\documentclass[12pt]{article}

\usepackage[utf8]{inputenc}
\usepackage{amsmath,amsfonts}
\usepackage{graphicx,epstopdf}
\usepackage{color,hyperref}  %  <-- to prevent option clash due to arXiv auto importing hyperref
\usepackage{subcaption}
\usepackage{xfrac}

\usepackage{pgfplots}
\pgfplotsset{compat=newest}

\newcommand{\Q}{\mathbf Q}
\newcommand{\Z}{\mathbf Z}
\newcommand{\reals}{\mathbf R}

 \newcommand{\R}{\mathbf R}

\newcommand{\pgl}{\mathrm{PGL}_2(\Z)}

\newcommand\nt[1]{\textcolor{red}{{\small { #1}}}}
\renewcommand\nt[1]{}  %  <--  disable notes

\newcommand{\Jimm}{\mathbf J}

\newtheorem{theorem}{Theorem}

\newtheorem{lemma}[theorem]{Lemma}

\title{Testing the transcendence conjectures of a modular involution of the real line and its continued fraction statistics}
\author{Hakan Ayral, A. Muhammed Uludağ}
\date{}

\newcommand{\sherh}[1]{\fboxsep=0pt\setlength{\fboxrule}{1pt}
\begin{center}
   \fbox{\colorbox{green}{
         \begin{minipage}[t]{13cm}
            #1
         \end{minipage}
      }
   }
\end{center}}
\newcommand{\sherhh}[1]{\fboxsep=0pt\setlength{\fboxrule}{1pt}
\begin{center}
   \fbox{\colorbox{yellow}{
         \begin{minipage}[t]{13cm}
            #1
         \end{minipage}
      }
   }
\end{center}}

\newcommand{\sherhhh}[1]{\fboxsep=0pt\setlength{\fboxrule}{1pt}
\begin{center}
   \fbox{\colorbox{red}{
         \begin{minipage}[t]{13cm}
            #1
         \end{minipage}
      }
   }
\end{center}}

\renewcommand{\sherh}[1]{}\renewcommand{\sherhh}[1]{}\renewcommand{\sherhhh}[1]{}%\renewcommand{\nt}[1]{}\renewcommand{\unut}[1]{}

\begin{document}

\maketitle

\begin{abstract}
We study the values of the recently introduced involution $\Jimm$ (jimm) of the real line, which is equivariant with the action of the group PGL(2,Z).
We test our conjecture that this involution sends algebraic numbers of degree at least three to transcendental values.  We also deduce some theoretical results concerning the continued fraction statistics of the generic values of this involution and compare them with the experimental results.
\end{abstract}

\section{Introduction}
Every irrational real number admits a unique simple continued fraction representation 
$$
x=n_0+\cfrac{1}{n_1+\cfrac{1}{n_2+\cfrac{1}{\dots}}} 
\quad (n_0\in \Z, n_1, n_2\dots \in \Z_{>0})
$$
denoted shortly as $x=[n_0,n_1,n_2,\dots]$. The numbers $n_0, n_1, \dots $ are called the partial quotients of $x$. If $x$ is a real quadratic irrational, i.e. if $x=a+\sqrt{b}$ with $a,b\in \Q$ and $b>0$ is a non-square, then $x$ is known to have an eventually periodic continued fraction representation. In contrast with this, not much is known about the continued fraction representations of other real algebraic numbers.

Denote by $\overline{\Q}$ the field of algebraic numbers, i.e. the set of roots of polynomials with coefficients in $\Q$. If $x\in \overline{\Q}\cap \R$ is not quadratic, numerical evidence  suggest that this representation should behave like the expansion of a ``normal" number (i.e. its partial quotients must obey the Gauss-Kuzmin statistics). In particular, the partial quotients averages are expected to tend to infinity. However, to our knowledge the answer to the much weaker question ``{are partial quotients of $x\in \overline{\Q}\cap \R$  unbounded if $n$ is not quadratic?}'' is currently not known \cite{adam1}. 

We have recently introduced and studied (\cite{dynamicaljimm}, \cite{subtlesymmetry}, \cite{shortjimm}) a certain 
continued fraction transformation $$\Jimm:\R-\Q \to \R,$$ which is involutive (e.g. $\Jimm\circ \Jimm =Id$) and  sends 
{normal numbers}   to ``{sparse numbers}" (i.e. numbers with partial quotients equal to 1 with frequency 1).  To define $\Jimm(x)$ for irrational $x=[0,n_1,n_2,\dots]$ in the unit interval $[0,1]$, assume first $n_1, n_2, \dots \geq 2$. Then 
\begin{equation}\label{def}
\Jimm([0,n_1,n_2,\dots]):=[0,1_{n_1-1},2,1_{n_2-2},2,1_{n_3-2},\dots],
\end{equation}
where $1_k$ denotes the sequence ${1,1,\dots, 1}$ of length $k$.
To extend this definition for $n_1, n_2, \dots \geq 1$, we eliminate the $1_{-1}$'s emerging in (\ref{def}) by  the rule $[\dots m, 1_{-1},n,\dots]\to [\dots m+n-1,\dots]$ and $1_0$'s by the rule 
$[\dots m, 1_{0},n,\dots]\to [\dots m,n,\dots]$. These rules are applied once at a time. See the next section for more details on the definition and properties of $\Jimm$.

Recall that the degree of an algebraic number $x\in\overline{\Q}$ is defined to be the degree of the polynomial of minimal degree satisfied by $x$. Our aim here is to experimentally confirm the following conjecture: 

\bigskip\noindent\textbf {Transcendence conjecture:}
If $x$ is a real algebraic number of degree $>2$, then $\Jimm(x)$ is transcendental.

\bigskip
Here is the chain of reasoning which led us to this conjecture:
\begin{align}
x\in \overline{\Q}, \deg(x)>2& \implies  
x \mbox{ is typical (belief) } \label{R1}\\
x \mbox{ is typical} & \implies  \lim_{k\to \infty} \frac{n_1+\dots +n_k}{k} =\infty \label{R2}\\
\lim_{k\to \infty} \frac{n_1+\dots +n_k}{k} =\infty & \implies y=\Jimm(x)  \mbox{ is sparse}\label{R3}\\
y  \mbox{ is sparse} & \implies y \mbox{ is not algebraic (belief)} \label{R4}
\end{align}

In this reasoning, the statement (\ref{R1}) is based on the numerical evidence mentioned in the introduction and it is widely believed to be true.
The statement (\ref{R2}) is Khinchine's theorem, (\ref{R3}) is an easy observation (see Lemma 1 below) and finally (\ref{R4}) is a contrapositive instance of (\ref{R1}).

 Although there are some recent results in the literature, pertaining to the transcendence of sparse continued fractions, 
(see the works of Adamczewski \cite{adam1} and Bugeaud \cite{bugeaud}) we don't know how this conjecture can be proven. There is also a much bolder version of the transcendence conjecture. Recall that the $\pgl$ is the group of invertible linear fractional transformations of dimension 2, i.e.
$$
\pgl:=\left\{\frac{ax+b}{cx+d}\,:\, a,b,c,d\in \Z, \quad  ad-bc=\pm 1\right\}
$$
where the group operation is the functional composition. $\pgl$ acts naturally on 
the extended real line $\R\cup \{\infty\}$.

\bigskip\noindent\textbf {Strong transcendence conjecture:}
In addition to the transcendence conjecture, any set of algebraically related numbers in the set
$$S:=\{\Jimm(x)\,: \,x\in \overline{\Q}, \, \deg(x)>2\}$$ are in the same $\pgl$-orbit. \nt{need a better formulation}

\bigskip
As an example, $\Jimm(x)$ and $\Jimm(\frac{ax+b}{cx+d})$ are (provably) algebraically dependent for any $x\in \R$ and $\frac{ax+b}{cx+d}\in \pgl$, whereas if $x\in \overline{\Q}$ is non-quadratic, then  $\Jimm(x)$, $\Jimm(x^2)$ and $\Jimm(2x)$ are (conjecturally) not. A challenge might be to find some {\it real} number $x$ which is not rational nor a quadratic irrational, and such that $\Jimm(x)$, $\Jimm(2x)$ and $\Jimm(x^2)$ are algebraically dependent.

The next section of the paper is devoted to the involution $\Jimm$. Section 3 gives a theoretical study of frequencies of partial quotients of $\Jimm(x)$ for general $x$ and compare them with the experimentally found frequencies for  $\Jimm$-values for algebraic $x$. We also consider some transcendental $x$ such as the number $\pi$. 

% HHH
Our experiment methodology relies on exploration of algebraic numbers at the proximity of conjectured to be transcendent numbers through the lattice reduction algorithm PSLQ. We also investigate whether involution $\Jimm$ preserves the algebraic dependencies by searching for integer relations between the images of algebraically dependent number pairs, once again using PSLQ algorithm. Numeric results presented on the following sections are computed with Python using SymPy \cite{sympy} for symbolic computations,  mpmath \cite{mpmath} for arbitrary precision floating point arithmetic which in turn benefits from GNU MP \cite{gnump} through Gmpy2 wrapper. It must be stressed that the numerical transcendence tests of this paper are with high confidence though not with certainty.

\section{The involution $\Jimm$}
The involution $\Jimm$ originates from the  outer automorphism group of $\pgl$ and can be viewed as an automorphism of the Stern-Brocot tree of continued fractions. These descriptions immediately shows that $\Jimm$ is involutive. It has some very peculiar analytic, arithmetic and dynamical properties.

Here we give an overview of its definition and some of its properties. For details we refer to \cite{subtlesymmetry}, \cite{dynamicaljimm} and \cite{shortjimm}.

\subsection{Definition of the involution}
\begin{figure}[h!]
\centering
\noindent{\includegraphics[scale=0.3,trim=0cm 4cm 0cm 3.5cm]{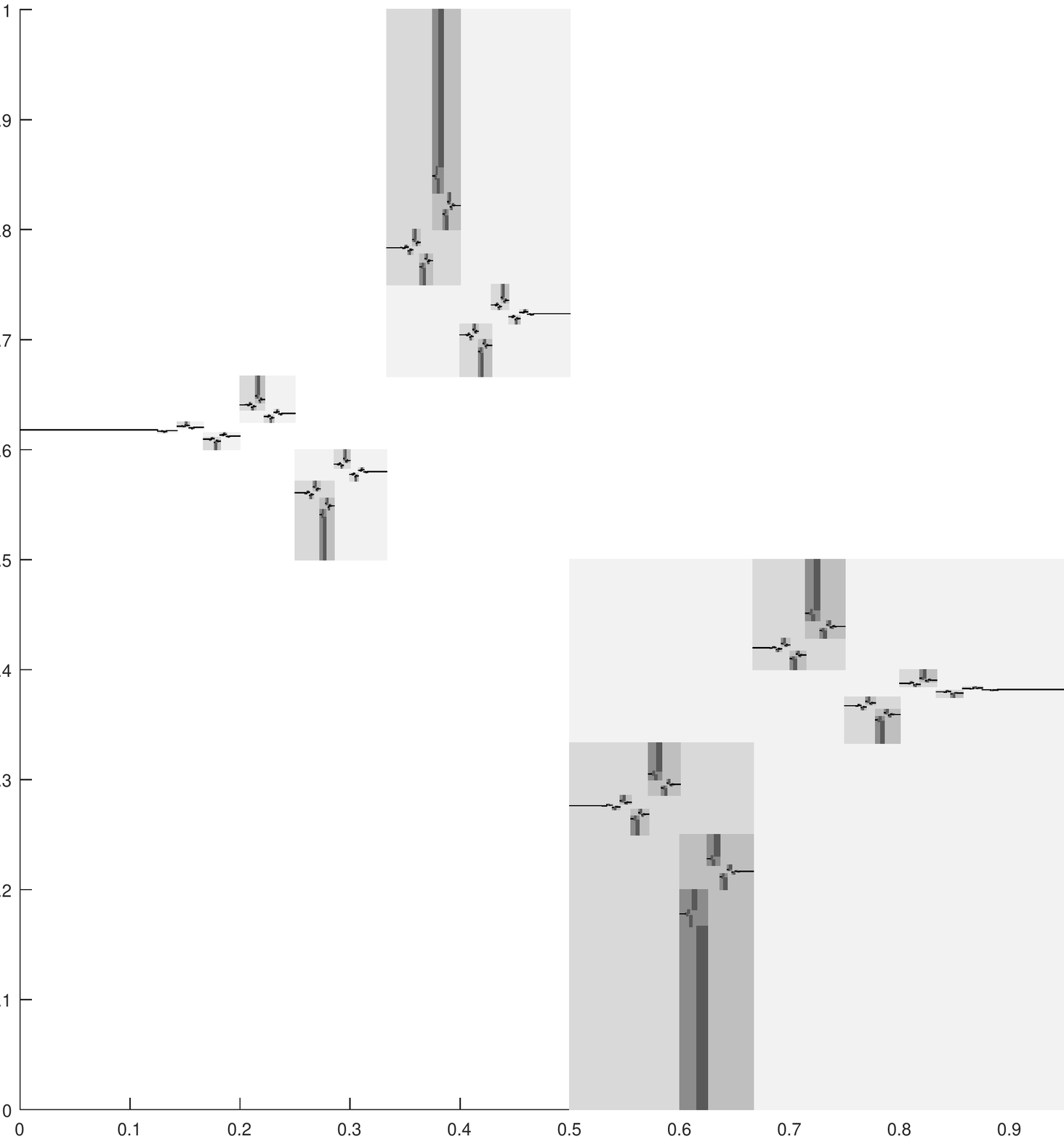}}
\caption{The plot of the involution $\Jimm$ on the unit interval.}

% \nt{1) mavi hatlar yok edilmeli 2) gradient yumusatılmalı; cok cabuk mutlak siyaha gidip kayboluyoruz. Mutlak siyah çizgi gibi görülmeli} 
\end{figure}

To illustrate the definition of $\Jimm$ given in the introduction, consider the following example:

\medskip\noindent
{\bf Example 1.} One has
\begin{align*}
\Jimm([1,
\textcolor{blue}{1},
\textcolor{red}{1},
\textcolor{green}{1},
\textcolor{magenta}{1},
13,\dots])=&
[ 1_{0},
\underbrace{2,
\textcolor{blue}{1_{-1}},
2}{},
\textcolor{red}{1_{-1}},
2,
\textcolor{green}{1_{-1}},
2,
\textcolor{magenta}{1_{-1}},
2,
1_{11}, \dots]\\
&\quad\,\,\, \quad\,\,[
\underbrace{3,\textcolor{red}{1_{-1}},
2}{},
\textcolor{green}{1_{-1}},
2,
\textcolor{magenta}{1_{-1}},
2,
1_{11}, \dots]\\
&\quad\,\, \quad\,\,\, \quad\,\,[
\underbrace{4,\textcolor{green}{1_{-1}},
2}{},
\textcolor{magenta}{1_{-1}},
2,
1_{11}, \dots]\\
&\quad\,\, \quad\,\,\, \quad\,\, \quad\,\,\, [ 
\underbrace{5,\textcolor{magenta}{1_{-1}},
2}{},
1_{11}, \dots]\\
&\quad\,\, \quad\,\,\,\,\, \quad\,\, \quad\,\, \quad\,\, [6,1_{11}, \dots]
\end{align*}

The next example should convince the reader that $\Jimm$ is involutive:

\medskip\noindent
{\bf Example 2.} One has
\begin{align*}
\Jimm([6,1_{11}, \dots])&=[1_5, 2, 1_{-1}, 2, \dots , 1_{-1},2, \dots]\\
&=[1,1,1,1,1,13,\dots]
\end{align*}

\paragraph{$\Jimm$ on rationals.} 
Since every nonzero rational number admits two simple continued fraction representations (one ending with a $1$ and the other not), and since the defining rule of $\Jimm$ produces two different values when applied to these representations, our involution is not well-defined on $\Q$.
In fact, there is a way to extend $\Jimm$ to the set of rationals as well, see \cite{subtlesymmetry}. 
In Appendix I, we provide a maple code which evaluates $\Jimm$ at a given rational.

\paragraph{$\Jimm$ on noble numbers.} 
A number is said to be {\it noble} if its continued fraction terminates with $1_\infty$.
If $x$ is noble, then $\Jimm(x)$ is rational. In fact
$$
\Jimm([n_0, \dots, n_k,1_\infty])=[1_{n_0-1},2,\dots,2,1_{n_k-2}]=
\Jimm([n_0, \dots, n_k-2,2,1_\infty]),
$$
so that $\Jimm$ is 2-1 on the set of noble numbers. 
%$\mathcal N$ 

\paragraph{$\Jimm$ on quadratic irrationals.} 
It is well known that ultimately periodic continued fractions are precisely real quadratic irrationals.
Since by definition $\Jimm$ preserves the periodicity of continued fractions, 
$\Jimm$ sends quadratic irrationals to quadratic irrationals (setwise). As an example,
 $\sqrt{2}=[1,\overline{2}]\implies\Jimm(\sqrt{2})=[\overline{2}]=1+\sqrt{2}$.
In general the situation is not so simple; we give a list of $\Jimm$-transforms of some quadratic surds in Appendix II below. This action respects the Galois conjugation, i.e. 
\begin{eqnarray}
\Jimm(a+\sqrt{b})=c+\sqrt{d}\implies \Jimm(a-\sqrt{b})=c-\sqrt{d}\\\Jimm(a+\sqrt{b})=c-\sqrt{d}\implies \Jimm(a-\sqrt{b})=c+\sqrt{d},
\end{eqnarray}
where $a,b,c,d \in \Q$ with $c,d>0$ being non-squares. This fact, together with the functional equations (see Section~\ref{property} below) below implies that 
if $x=\sqrt{q}$ where $q\in \Q$ is a positive non-square, with
$\Jimm(x)=a\pm\sqrt{b}$, then $\mathrm{Norm}(x):=a^2-b=-1$.
There is a host of such correspondences of quadratic irrationals under 
$\Jimm$, for example, 
$\mathrm{Norm}(x)=1\iff {\mathrm{Norm}(\Jimm (x))=1}$. See \cite{shortjimm} for details.

\paragraph{$\Jimm$ on other numbers.} \label{property}
To finish, let us give the $\Jimm$-transform of two numbers that we  used in our experiments:
\begin{eqnarray*}
\Jimm(\sqrt[3]{2})=
\Jimm([1; 3, 1, 5, 1, 1, 4, 1, 1, 8, 1, 14, 1, 10, 2, 1, 4, 12, 2, 3, 2, 1, \dots])\\
=[2,1,3,1,1,1,4,1,1,4,1_6,3,1_{12},3,1_8,2,3,1,1,2,1_{10},2,2,1,2,\dots ]
\\
=2.784731558662723\dots,\\
\Jimm(\pi)=
\Jimm([3, 7, 15, 1, 292, 1, 1, 1, 2, 1, 3,  \dots])=
[1_2, 2, 1_5, 2, 1_{13}, 3, 1_{290}, 5, 3,  \dots]\\
=1.723770792548027\dots.
\end{eqnarray*}

\subsection{Further properties of $\Jimm$}
\paragraph{Functional equations.}
We have defined $\Jimm$ on the unit interval. We may extend it to $\R$ via the equation $\Jimm(-x)=-1/\Jimm(x)$. This extension satisfies the functional equations (see \cite{shortjimm})
\begin{equation}\label{fes}
\Jimm(1/x)=1/\Jimm(x),\quad \Jimm(1-x)=1-\Jimm(x),\quad \Jimm(-x)=-1/\Jimm(x).
\end{equation}
In fact, these equations characterize $\Jimm$. 
Given the continued fraction representation of $x$, one can use these functional equations to compute $\Jimm(x)$. 
\paragraph{Harmonic numbers and Beatty partitions.}
Using the second and the third equations one obtains
$$
\frac1x+\frac1y=1 \iff \frac{1}{\Jimm(x)}+\frac{1}{\Jimm(y)}=1;
$$
so $\Jimm$ preserves harmonic pairs of real numbers. This implies that $\Jimm$ acts on the Beatty partitions (see \cite{stolarsky1976beatty}) of the set of natural numbers.
\paragraph{Modularity.}
Another consequence of the functional equations is that, if the continued fractions of $x$ and $y$ have the same tail, then the same is true for
 $\Jimm(x)$ and  $\Jimm(y)$. This shows that $\Jimm$ sends $\pgl$-orbits to $\pgl$-orbits, i.e. it defines an involution of $\R/\pgl$, the ``{moduli space of degenerate rank-2 lattices}". As such, we may consider it as a kind of ``{modular form}". It is easy to see that $\Jimm$ is continuous on $\R\setminus \Q$ and with jumps at rationals;
 and we were able to prove that $\Jimm$ is differentiable almost everywhere with a derivative vanishing almost everywhere \cite{derivativeofjimm}.
\paragraph{Dynamics.}
Recall that the celebrated Gauss continued fraction map ${\mathbb T}_G: [0,1]\to  [0,1]$
is the map which forgets the first partial quotient:
\begin{equation}
{\mathbb T}_G: x=[0,n_{1}, n_{2},n_{3},\dots]=[0,n_{2}, n_{3},n_{4},\dots] 
\end{equation}
Our involution $\Jimm$ conjugates the Gauss continued fraction map to the so-called Fibonacci map ${\mathbb T}_{F}:[0,1]\to  [0,1]$, defined as
\begin{equation}
\Jimm \circ {\mathbb T}_G \circ \Jimm:{\mathbb T}_{F}: 
[0,1_k,n_{k+1},n_{k+2}, \dots] \to 
[0,n_{k+1}-1,n_{k+2}, \dots]
\end{equation}
where it is assumed that $n_{k+1}>1$ and $0\leq k <\infty$. Dynamical properties of ${\mathbb T}_{F}$ and ${\mathbb T}_{G}$ are tightly related, (see \cite{isola2014continued} and \cite{dynamicaljimm}). For example the eigenfunctions of their transfer operators (see \cite{mayer2012iv}) satisfies the same three-term functional equation studied in \cite{lewis2001period}.

Our hope is that, due to these rich properties of the involution $\Jimm$,  especially the functional equations (\ref{fes}) and its effect on quadratic irrationals, it might be possible to infer the transcendence of $\Jimm(x)$ directly from the knowledge of algebraicity of $x$; by-passing the ``{beliefs}" in our chain of reasonings which led us to the transcendence conjecture.

\section{The transcendence conjectures}
Let us explain the theoretical basis for the conjecture. If $X\in [0,1]$ is a uniformly distributed random variable, then by Gauss-Kuzmin's theorem \cite{khinchinebook}, the frequency of an integer $k>0$ among the partial quotients of $X$ equals almost surely
$$
p(k)=\frac{1}{\log 2}\log \left(1+{1\over k(k+2)}\right) .
$$
\sherh{{\bf Exercise}. Given two probability distributions $p(n)$, $q(n)$ on $\Z_{>0}$, find $x=[0, n_1, n_2, \dots]\in [0,1]$ with $\Jimm x=[0, m_1, m_2, \dots]$ such that 
$$\lim_{k\to\infty} { \sharp \{1\leq i\leq k\, :\, n_i=n\} \over k}=p(n)$$
\mbox{ and }
$$\lim_{k\to \infty} { \sharp \{1\leq i\leq k\, :\, m_i=m\} \over k}=q(m)$$}

\noindent
Moreover, the arithmetic mean of its partial quotients tends almost surely to infinity (see \cite{khinchinebook}), i.e. if 
$X =[0, n_1, n_2, \dots]$ then 
\begin{equation}
\lim_{k\to \infty} {n_1 +\cdots +n_k \over k} = \infty\qquad\hbox{(a.s.)}
\end{equation}
In other words, the set of numbers in the unit interval such that the above limit exists and is infinite, is of full Lebesgue measure. Denote this set by $\mathcal A$. Since the first $k$ partial quotients of $X$ give rise to at most $n_1 +\cdots +n_k -k$ partial  quotients of $\Jimm(X)$ and 
at least $n_1 +\cdots +n_k -2k$ of these are 1's, one has 
$$
{n_1 +\cdots +n_k -k \over  n_1 +\cdots +n_k -2k} \to 
{\frac{n_1 +\cdots + n_k}{k} -1 \over  \frac{n_1 +\cdots + n_k}{k} -2} \to 1
$$
This proves that $\Jimm(x)$ is almost surely `{sparse}' in the following sense:
\begin{lemma}\label{density}
The frequency of $1$'s among the partial quotients of $\Jimm(X)$ equals 1 a.s.. In particular the partial quotient averages of $\Jimm(X)$ tend to 1 a.s. and $\Jimm(\mathcal A)$ is a set of zero measure.
\end{lemma}
Note that $x$ and $\Jimm(x)$ can be simultaneously sparse, consider e.g. $x=[0,1_{2}, 2^2, 1_{2^3}, 2^4, 1_{2^5}\dots]$.

 Below is a partial quotient statistics for 
the numbers $\Jimm(\pi)$ and $\Jimm(\sqrt[3]2)$. (see Table\ref{tab:table-jpi-vs-jcbrt2})

\begin{table}[h!]
\centering
\begin{subfigure}[t]{0.31\textwidth}
	\begin{tabular}{ l l}
	partial \\ quotient  & frequency\\
	\hline

1  & 95.160  \\ 2  & 2.636  \\ 3  & 1.418  \\ 4  & 0.471  \\ 5  & 0.186  \\ 6  & 0.078  \\ 7  & 0.033  \\ 8  & 0.009  \\ 9  & 0.004  \\ 11  & 0.002  \\ 10  & 0.001  \\ 13  & 0.001  \\

	\end{tabular}
	\caption{Statistics of $\Jimm(\pi)$.}

\end{subfigure}
\hspace{1cm}
\begin{subfigure}[t]{0.31\textwidth}
	\begin{tabular}{ l l }
	partial \\ quotient  & frequency\\
	\hline

1  & 94.761  \\ 2  & 2.891  \\ 3  & 1.535  \\ 4  & 0.476  \\ 5  & 0.207  \\ 6  & 0.073  \\ 7  & 0.034  \\ 8  & 0.013  \\ 9  & 0.004  \\ 11  & 0.001  \\ 10  & 0.000  \\ 13  & 0.000  \\

	\end{tabular}
	\caption{Statistics of $\Jimm(\sqrt[3]2)$.}

\end{subfigure}
\caption{Both statistics have been made for the continued fraction expansions of $\pi$  of length $10^5$ and of $\sqrt[3]2$ of length $2 \times 10^4$. 
For mentioned expansion sizes, the length of $\Jimm(\pi)$ is $> 1.4 \times 10^6$ terms whereas that of $\Jimm(\sqrt[3]2)$ is $>2.2 \times 10^5$.}  %  is $1,426,400$ terms whereas that of $\Jimm(\sqrt[3]2)$ is $223,642$.}
\label{tab:table-jpi-vs-jcbrt2}
\end{table}

The proof of the Lemma applies to any number $x$ with partial quotient averages tending to infinity.
Since for  $x$ algebraic of degree $>2$,  it is widely believed that this is the case, we see that the average partial quotient of 
$\Jimm(x)$ is very likely to tend to 1 for such $x$. Since this is far from being unbounded, it is natural to believe that the transcendence conjecture is true.

\section{Testing transcendence by searching algebraic numbers on neighborhood of $\Jimm(\sqrt[3]2)$}

In this section we provide some lower bounds, in terms of degree and coefficient size, for a possible minimal polynomial which may have the image of $\sqrt[3]2$ under $\Jimm$ as a root. We picked $\sqrt[3]2$ as the simplest representative for algebraics of degree $\ge 3$ which our conjecture deals with.

For computations we start with a continued fraction approximation of $\sqrt[3]2$, truncated at $2\times 10^4$ terms. After applying the involution we obtain an approximation to $\Jimm(\sqrt[3]2)$ with $223642$ partial quotients. The rational number that correspond to this c.f. expansion has denominator $> 10^{49552}$, hence our representation provide up to $49552$ decimal digits for the image under $\Jimm$.

Here the sparsity mentioned on previous section manifest itself as a lack of information content on numerical computations. Notice that each c.f. term of $\Jimm(\sqrt[3]2)$ provide $\frac{\ln(10) \times 49552}{\ln(2) \times 223642} \approx 0.736$ bits per partial quotient; this is much lower than the entropy of Gauss-Kuzmin distribution which is approximately $3.432$ bits \cite{gkentropy}. On the other hand partial quotients of ``{typical}" numbers are expected to obey Gauss-Kuzmin distribution, which implies having a specific information density different than what we observe on our example.

In order to gather some bounds on degree and coefficient size, we searched the polynomials with degree up to 32 and coefficients up to $10^{100}$ in absolute value, using the PSLQ algorithm which performs integer relation search by lattice reduction\cite{pslq1,pslq2}. For numeric computation we setup a working precision of $10^4$ decimal digits (i.e. $\epsilon = 10^{-10000}$) to keep the computation time manageable, even though our approximation of $\Jimm(\sqrt[3]2)$ is correct up to five times more decimal digits.

An integer relation is an equation of the form $a_1x_1 + a_2x_2 + \cdots + a_nx_n = 0$ with $a_i \in \Z$ and $x_i \in \reals$, and an integer relation algorithm is an algorithm which searches $a_i$ values that satisfies a such equation for a given sequence of $x_i$ values.

The well known way to use any integer relation algorithm to search for a polynomial up to a given degree and having a specific value as its root is as follows: we build a list of successive powers of the value, such as $[x^0,x^1,x^2,x^3,\dots,x^n]$ where $x$ is the desired root specified with enough precision, and $n$ is the maximum degree that we want to limit our search within. Setting the $x_i$ values as the powers of a single $x$ value effectively reduces the left side of the previous equation to evaluation of all polynomials (bounded by degree and coefficients) at point $x$.

Given a such list along with two values, one limiting the size of allowed coefficients in absolute value (say $c_{max}$), and the other representing the acceptable difference from zero as the termination condition of the algorithm (say $\epsilon$), PSLQ algorithm search smallest integer coefficients $a_i$ satisfying $|a_1 x_1 + a_2 x_2 + ... + a_n x_n| < \epsilon$.

To search all polynomials of degree 32 and below we successively invoke the algorithm with the set of input variables $\{x^n : 0 \le n \le m \}$ where $m$ is increased by one at each iteration up to the value 32. For each degree the algorithm returns the coefficients of the polynomial within the bounds that attain the smallest absolute value.

Take note that, just like rational approximation, with integer relations algorithms it is always possible to find a better algebraic approximation to any irrational either by increasing the maximum allowed degree for the polynomial, or by allowing larger coefficients.

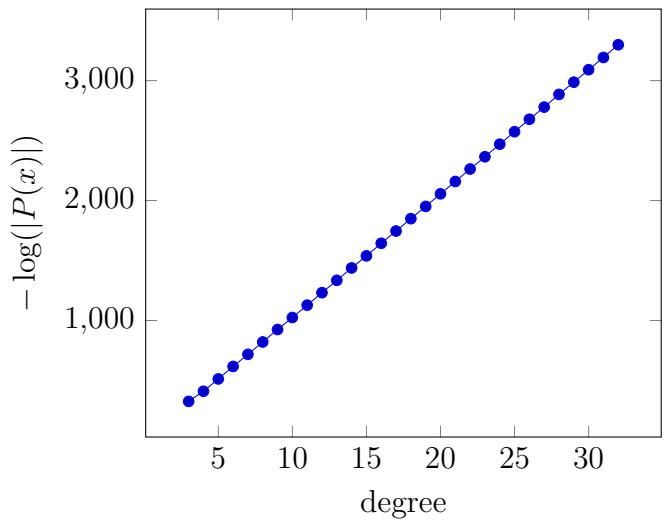
\begin{figure}
\begin{center}
\hspace{-1cm}
\begin{tikzpicture}
\begin{axis} [
xlabel={degree},
ylabel={$-\log(|P(x)|)$},
]
\addplot table [x=degree, y=mlog, col sep=tab] {data.txt};
\end{axis}

\end{tikzpicture}
\end{center}
\caption{Residual difference from zero}
\label{fig:errors}
 \end{figure}

In agreement with our transcendence conjecture, PSLQ algorithm did not find any polynomial having $\Jimm(\sqrt[3]2)$ as a root when bounded with degree up to 32 and coefficients up to $\pm10^{100}$; therefore for each degree we recorded the smallest absolute value attained by any polynomial before the algorithm gave up, along with the reported norm bound of the polynomials for which the algorithm guarantees no value closer to zero can be attained with. As expected, each time the reported bound on $L^\infty$ norm was very close to the coefficient upper bound of $10^{100}$ that we provided. For constant coefficient bound, the exponent of distance to zero exhibits a fairly linear relationship to the degree of polynomial (see \autoref{fig:errors}), this confirms that our computations obey to the relationship given by Bailey in \cite{bailey-pslq}, where he states that in order to recover an integer relation the required precision in number of digits is $n\log(G)$ where $n$ is the degree bound and $G$ is the bound on coefficients. This precision requirement is expected both for input values and for precision of intermediate numerical computations; as previously stated both the number of partial quotients we use to approximate $\Jimm(\sqrt[3]2)$ and our working precision for arithmetic operations are well above Bailey's minimum precision requirements, but the recovered relation can only approach to zero as much as predicted by Bailey's relationship due to constraints on coefficient size.

% article on pslq error

It must be noted that any finite representation of an irrational is eventually an approximation to the idealised number, unless the representation is an algorithmic one; thus the above mentioned finite continued fraction representation is effectively a rational approximation. On the other hand, amount of accuracy listed above is more than enough for the search space which is practically limited by the computational power. Also as the final polynomial we obtained by this procedure has 32 coefficients with 100 digits each, it is not practical to include it here.

%\eject
\section{Statistics}
Let $X\in [0,1]$ be a uniformly distributed random variable.
Lemma~\ref{density} implies that the density of any $k>1$ among the partial quotients of $\Jimm(X)$ is a.s. zero. 
What if we ignore the 1's in the partial quotients? 
To be more precise, define the ``{collapse of continued fraction map}"
as 
$$
{\mathsf P}: [0,n_1,n_2,\dots, ]\to [0,n_1-1,n_2-1,\dots],
$$
where the vanishing partial quotients are simply ignored. For example, 
$$
{\mathsf P}([0,2,4,1,1,1,9,5,\dots])=[0,1,3,8,4,\dots].
$$
(We don't care about the collapse of the continued fractions terminating with $1_\infty$ as these are countable in number.)
Then we can determine the partial quotient statistics
of ${\mathsf P}\Jimm(X)$, as follows. It is known that the frequency of $1_i$ in the continued fraction expansion of $X$
equals (see~\cite{hakami}) 
$$
k(i)=\frac{(-1)^i}{\log 2}\log\left(1+\frac{(-1)^i}{F_{i+2}^2}\right)
$$
This counting involves a certain repetitiveness in that, if $j\leq i$, then the string $\dots n, 1_i, m, \dots$ $(n,m>1)$ in the continued fraction contributes $k-l+1$ to the census. Denoting by $m(i)$ the frequency of the string $1_i$ occur, but not as a substring of a longer string of 1's, we have
$$
\left[\begin{matrix}
k_1\\
k_2\\
k_3\\
\dots\\
\end{matrix}
\right]
=
\left[\begin{matrix}
1& 2& 3&\dots &n\\
0&1&2&\dots&n-1\\
0&0&1&\dots&n-2\\
&&&\dots&\\
\end{matrix}
\right]
\left[\begin{matrix}
m_1\\
m_2\\
m_3\\
\dots\\
\end{matrix}
\right]
$$
$$
\implies
\left[\begin{matrix}
m_1\\
m_2\\
m_3\\
\dots\\
\end{matrix}
\right]
=
\left[\begin{matrix}
1& -2& 1&0&0&\dots &0&0\\
0&1&-2&1&0&\dots&0&0\\
0&0&1&-2&1&\dots&0&0\\
&&&&&\dots&\\
\end{matrix}
\right]
\left[\begin{matrix}
k_1\\
k_2\\
k_3\\
\dots\\
\end{matrix}
\right]
$$
Hence, we get $m_i=k_i-2k_{i+1}+k_{i+2}$. Note that
$$
\sum_{i=0}^\infty m(i)=0.584962500,\quad   \sum_{i=1}^\infty  im(i)=p(1)=0.41503749927,
$$
and that $m(0)$ is the frequency of the string $\dots, n,m, \dots$ with $n,m>1$. This string transforms under $\Jimm$ to the string $1_{n-1},2,1_{m-1}$, which after the collapse map yields a partial quotient 1 in ${\mathsf P}\Jimm(X)$.
The desired frequency of $i$ among the partial quotients
of ${\mathsf P}\Jimm(X)$ is the number $u(i)$ given by 
$$
u(i)=\frac{m(i-1)}{\sum_{j=0}^\infty m(i)}
$$  
Above we tabulate these theoretical values of the frequencies \autoref{tab:theoretical}, followed by the experimental values obtained by computation for $\mathsf P\Jimm(\pi)$ and $\mathsf P\Jimm(\sqrt[3]2)$ \autoref{tab:frequencies}.

\begin{table}[h!]
$$
\begin{array}{l|llll}
i&k(i)&m(i)&m(i)/\scriptsize{\sum }m_i&p(i) \mbox{\small (Gauss-Kuzmin)}\\\hline
0&1.0&                                             0.3219280948&0.5503397134&\\
1&                          0.4150374993&0.1699250016&0.2904887087&0.4150374989\\
 2&                         0.1520030934&0.0565835283&0.0967301805&0.1699250015\\
    3&                      0.0588936890&0.0227200765&0.0388402273&0.09310940485\\
         4&                 0.0223678130&0.0085115001&0.0145505056&0.05889368952\\
              5&            0.0085620135&0.0032751312&0.0055988737&0.04064198510\\
               6&           0.0032677142&0.0012474677&0.0021325602&0.02974734293\\
                   7&       0.0012485461&0.0004770014&0.0008154393&0.02272007668\\
                        8&  0.0004768458&0.0001821241&0.0003113433&0.01792190800\\
                        9&  0.0001821469&0.0000695769&0.0001189425&0.01449956955\\
                        10&0.0000695722&0.0000265739&0.0000454285&0.01197264119
                        \end{array}
$$
\caption{Theoretical values for the expected frequencies}\label{tab:theoretical}
\end{table}

%\eject

%%%The results of two experiments we performed are given in \autoref{tab:frequencies} .

\begin{table}[h!]

\centering
\begin{subfigure}[t]{0.31\textwidth}
	\begin{tabular}{ l l}
	partial \\ quotient  & frequency\\
	\hline

1  & 54.891  \\ 2  & 29.250  \\ 3  & 9.703  \\ 4  & 3.854  \\ 5  & 1.391  \\ 6  & 0.541  \\ 7  & 0.229  \\ 8  & 0.088  \\ 9  & 0.023  \\ 10  & 0.018  \\ 11  & 0.003  \\ 12  & 0.003  \\ %13 & 0.000

	\end{tabular}
	\caption{Statistics of $\mathsf P\Jimm(\pi)$.}

\end{subfigure}
\hspace{1cm}
\begin{subfigure}[t]{0.31\textwidth}
	\begin{tabular}{ l l }
	partial \\ quotient  & frequency\\
	\hline

1  & 55.202  \\ 2  & 29.304  \\ 3  & 9.099  \\ 4  & 3.952  \\ 5  & 1.408  \\ 6  & 0.657  \\ 7  & 0.256  \\ 8  & 0.085 \\  9  & 0.008  \\ 10 & 0.025  \\ 11  & 0.000  \\ 12  & 0.000  \\ %13  & 0.000  \\

	\end{tabular}
	\caption{Statistics of $\mathsf P\Jimm(\sqrt[3]2)$.}

\end{subfigure}

\caption{Digit frequencies observed by numerical computation}\label{tab:frequencies}
\end{table}

%\nt{Hakan! this is for you.. Recipe. In fact we already did this. Take $x=^3\!\!\!\sqrt{2}$ or $x=\pi$. Find $\Jimm(x)$. Then find $\mathsf P\Jimm(x)$ and compute its partial quotient statistics. Hopefully they will match with the above table.}}  %  <--  DONE

If we consider $\mathsf P\Jimm$ as a kind of derivation, then it is possible to compute the statistics of higher derivatives, by using results of (\cite{fergy}) where  the Gauss-Kuzmin statistics for $n_k$ $(n,k=1,2,\dots)$ have been computed.

\subsection{Conjectures concerning algebraic operations}  %  ----------------------------------------  \subsection{Algebraic operations}

%%  ---  2*J(pi)  ,  2*J(cbrt(2))  ---
%\begin{table}[h!]
%\centering
%\begin{subfigure}[t]{0.31\textwidth}
%	\begin{tabular}{ l l}
%partial \\ 
%quotient  & frequency\\
%\hline
%1  & 7.796  \\ 
%2  & 3.290  \\ 
%3  & 4.158 \\ 
%4  & 81.280 \\  
%5  & 0.880  \\   
%6  & 1.025 \\       
%7  & 0.407     \\ 
%8  & 0.609 \\   
%9  & 0.138 \\  
%10  & 0.194 \\ 
%11  & 0.063  \\ 
%12  & 0.075  \\ 
%13 & 0.020
%\end{tabular}
%\caption{Statistics of $\mathsf 2 \times \Jimm(\pi)$.}
%\end{subfigure}
%\hspace{.5cm}
%\begin{subfigure}[t]{0.31\textwidth}
%	\begin{tabular}{ l l }
%partial \\ 
%quotient  & frequency\\
%\hline
%1 & 9.496 \\
%2 & 3.947 \\
%3 & 4.940 \\
%4 & 77.376 \\
%5 & 1.081 \\
%6 &   1.295 \\
%7 & 0.512 \\
%8 & 0.723  \\
%9 & 0.148 \\
%10 & 0.200 \\
%11 & 0.067 \\
%12 & 0.106 \\
%13 &  0.028
%\end{tabular}
%	\caption{Statistics of $\mathsf 2 \times \Jimm(\sqrt[3]2)$.}
%
%\end{subfigure}
%\caption{Digit frequencies observed by numerical computation}
%%%%%%   \nt{pls change the ordering}}  %  <--- DONE
%\end{table}
%
%
%
%
%
%\nt{(I have no idea about the theoretical frequencies in this situation) Hakan! We also did this. You may include these if you have time:
%Do the statistics of $n\times (\Jimm(x))$ $(N=2,3,\dots)$ and see if they obey some distribution.
%}

Recall that $\mathcal A\subset [0,1]$ is the set of real numbers whose partial quotient averages tend to infinity
and that it is of full measure. It includes the set of ``{normal}" numbers where by ``{normal}" we mean that $x$ and $y$ obeys all predictions of the Gauss-Kuzmin statistics.

\medskip\noindent\textbf {Conjecture 3.}
 If $x, y\in \mathcal N$, then $\Jimm(x)+\Jimm(y)$ and $\Jimm(x)\Jimm(y)$ are normal a.s.. 

\medskip
Note that if $x$ is normal then so is $1-x$, whereas $\Jimm(x)+\Jimm(1-x)=1$ is surely not normal.
Also note that the set $\Jimm(\mathcal A)$ and therefore the sets $\Jimm(\mathcal A)+\Jimm(\mathcal A)$ 
and  $\Jimm(\mathcal A)\times \Jimm(\mathcal A)$ are of zero measure.

\begin{figure}[h!]
\centering
\includegraphics[scale=0.5,trim={2cm 1cm 2cm 1.5cm},clip]{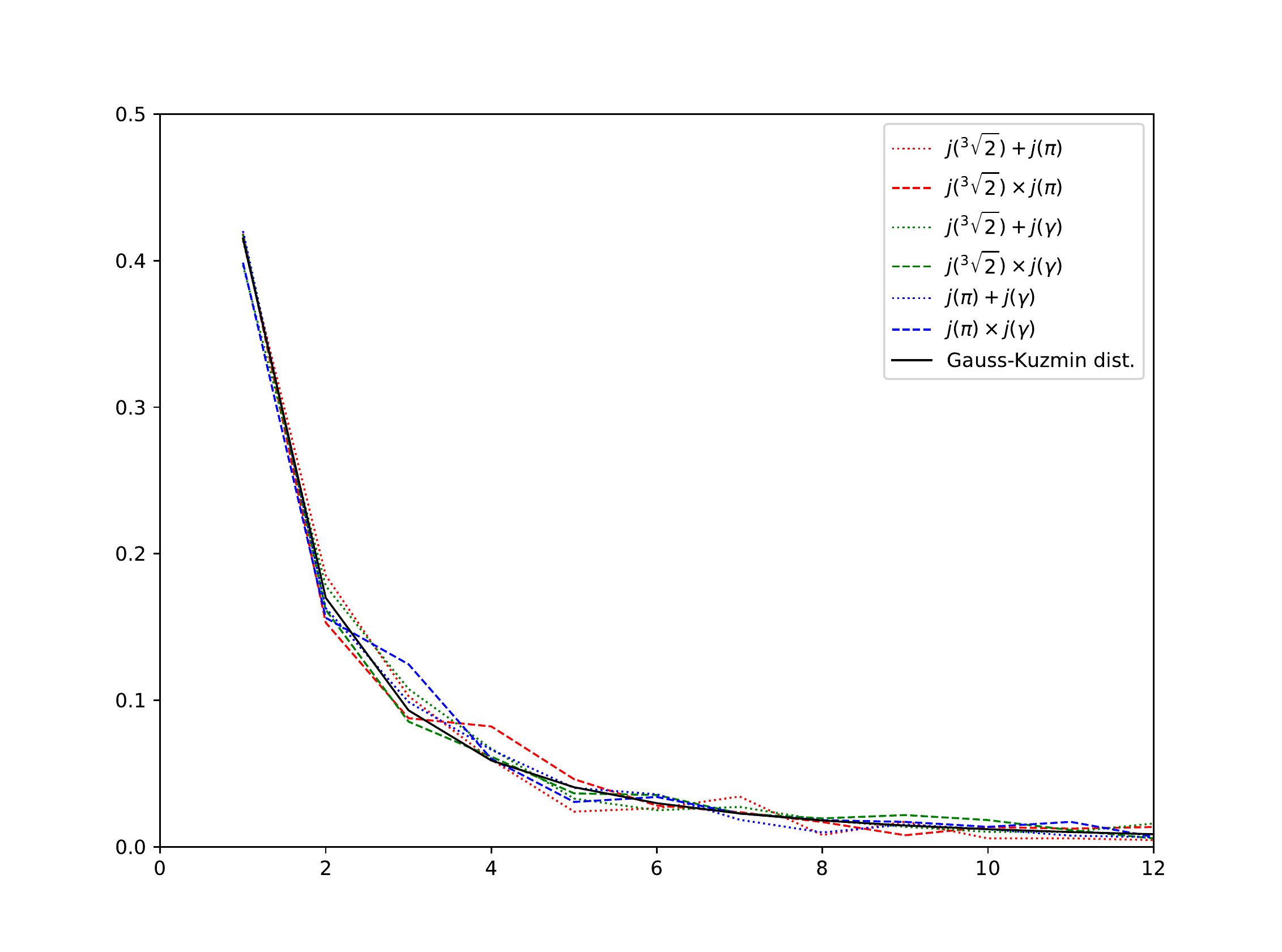} % generated with main6.py
\caption{Partial quotient frequencies of some presumably normal numbers under algebraic operations ($\gamma$ is Euler-Mascheroni constant)}
\label{fig:conj2}
\end{figure}

On \autoref{fig:conj2} the partial quotient frequencies for the sums and multiplications between $\Jimm(\pi),\Jimm(\sqrt[3]{2})$ and $\Jimm$-transformed Euler-Mascheroni constant $\Jimm(\gamma)$ are portrayed, along with the Gauss-Kuzmin distribution as a reference for the expected term distribution for normal numbers. A stem plot instead of the line plot would better suit for displaying term frequencies as these are only defined for integers; but as we want to show the overlap between multiple distributions, stem plot would end up with a too occluded graph, hence the use of line graph instead. It can be seen that even though the $\Jimm$ transformed numbers are `sparse' on themselves, any algebraic operation between them result with a term distribution which is characteristic of normal numbers.

\medskip\noindent\textbf {Conjecture 4.} $q\Jimm(X)$ obeys a certain law for every $q\in \Q$.

\begin{figure}[h!]
\centering
\includegraphics[width=\textwidth]{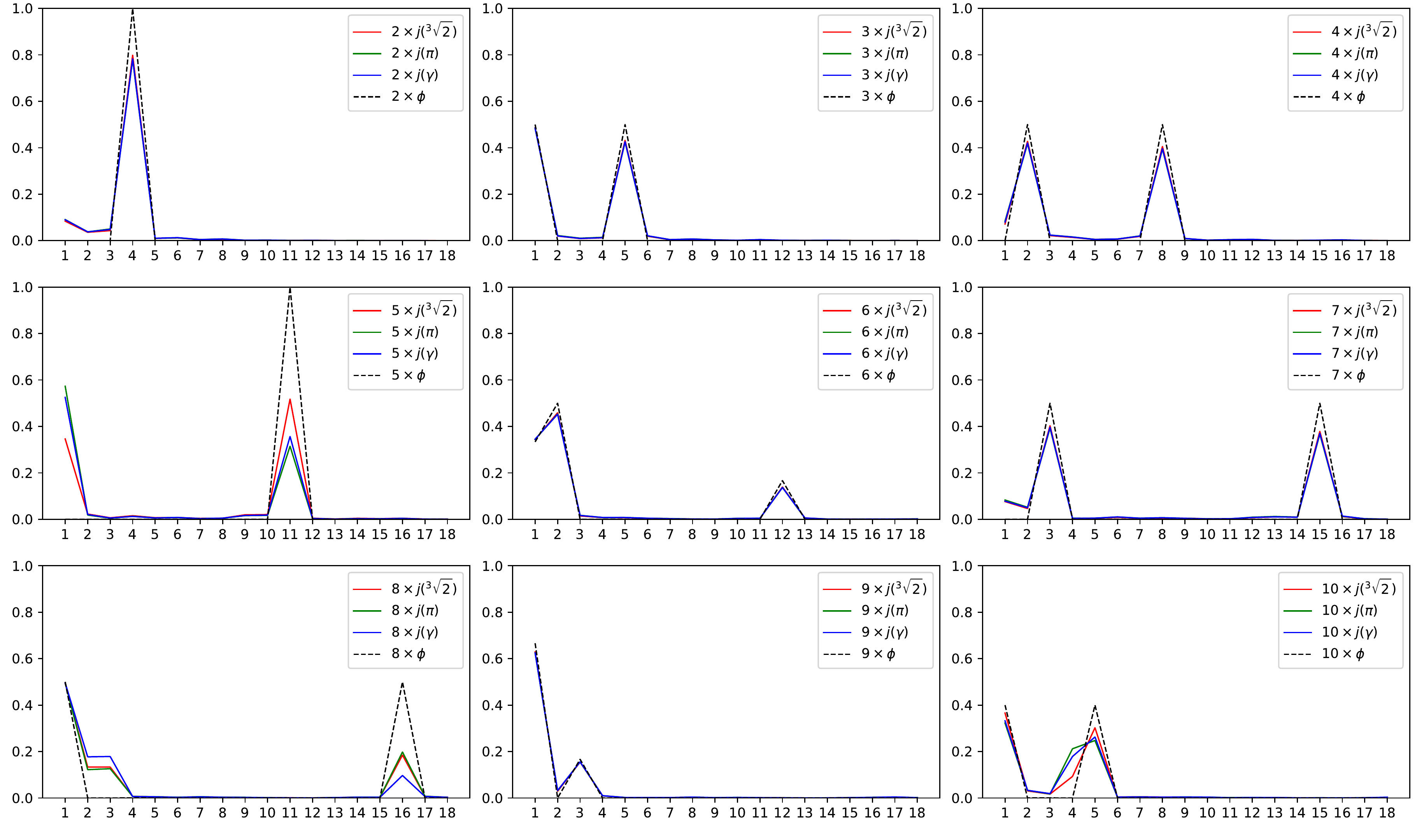} 
\caption{Partial quotient frequencies of some presumably normal numbers under algebraic operations ($\gamma$ is Euler-Mascheroni constant)}
\label{fig:conj-jtz}
\end{figure}

\begin{figure}[h!]
\centering
\includegraphics[width=\textwidth]{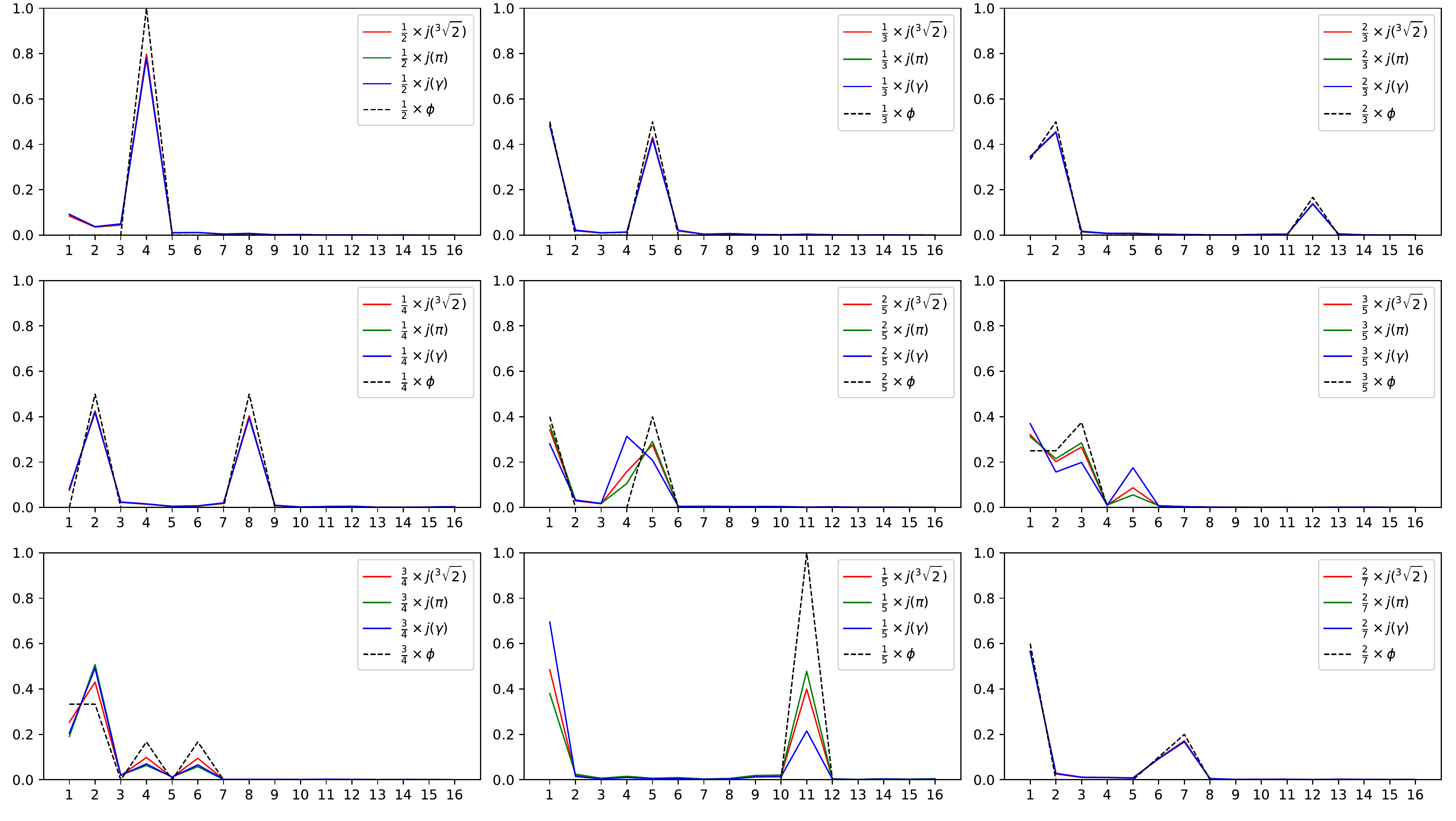} 
\caption{Partial quotient frequencies of some presumably normal numbers under algebraic operations ($\gamma$ is Euler-Mascheroni constant)}
\label{fig:conj3}
\end{figure}

\medskip\autoref{fig:conj-jtz} and \autoref{fig:conj3} display a few cases cases of the form $q\Jimm(X)$ with $q \in \Z$ and $q \in \Q$ respectively. In both figures the distributions of $\Jimm$ transformed (presumably) normal numbers and golden ratio are in agreement under multiplication with the same constant. Our study on the mechanism inducing this agreement is still in progress.

\nt{There is a possibility that there will be statistics for $q\Jimm(X)$, which may be different from the statistics for 
$q\Jimm(\Phi-1)$ (although we have a reason to believe that these will be close). We need to do several experiments with different numbers to see if they have the same statistics. 
And that will be the conjecture:   
I don't know if this will be true for all sparse numbers or shall we get something specific to $\Jimm$.}

\nt{\textbf{Conjecture 4.}It seems of interest to test the conjecture: If $x$ is sparse (this is more general then being just the jimm of some normal number), then the statistics of $P(x)$ will be same as the statistics of $P([1_\infty])$. In particular, this conjecture says that if $P(x)=x^2-x-1$ and $x$ is sparse, then so is $P(x)$. Perhaps this is the most interesting special case of this conjecture. For number crunching. Note that at this point we leave the context of $\Jimm$ and consider the much rough (and general) notion of sparse numbers. $\Jimm$ has a very delicate statistics whereas a sparse number has a very simple definition. To refine this notion, we may define $n$-sparse number to be a number with $\mathsf P^k x$ sparse for $1\leq k\leq n$, where $\mathsf P$ is the collapse operator.}

\nt{It might be possible to prove these results by carefully observing the Gosper algorithm.}

%%%%%
%ONE LAST EXPERIMENT
%$J(^3\sqrt{2}), J(^3\sqrt{4})$ are these numbers algebraically independent?

\subsection{Algebraic independence}

Another property to investigate is whether the algebraic dependence (or independence thereof) is a property that $\Jimm$ transformation preserves. In this regard we conducted a series of computations to search for possible relations between images of some algebraically related numbers under $\Jimm$. To achieve this we build a set consisting of the two images for which we investigate whether the dependence is preserved, some transformations of those (i.e. $\frac{1}{x},e^x,\log(x),\sqrt{x}, x^2$), some multiplicative combinations of those (i.e. $xy,\frac{x}{y},\frac{y}{x}$), along with transformations of combinations and combinations of transformations. This set serve as a dictionary of numbers which are algebraically related to one of the two $\Jimm$ images for which we investigate whether the original algebraic dependence before $\Jimm$ is carried in some form.

By running an integer relation algorithm on a such set we can recover more complex and non-linear relations than integer relations. We restricted the coefficients of the integer relation search by $\pm10^3$; once again it is always possible to approach arbitrarily close to any real number by allowing larger coefficients or by incorporating more complex transformations to the dictionary set.

\paragraph{Experiment 1.} In this experiment we search for a possible algebraic relation between $\Jimm(\sqrt[3]{2})$ and $\Jimm(\sqrt[3]{4})$; we progressively increase the allowed error tolerance in order to observe the complexity of proposed relation as a function of allowed error.

\begin{table}[h!]
\center
	\begin{tabular}{ l r }
\hline
	Tolerance & Relation Expression  \\ \hline \\

%$10^{-2}$  & \large $\alpha={\frac{-5}{3} + \beta}$  \\[0.4cm]

%$10^{-3}$  & \large $\alpha={\frac{14}{13} + \frac{5}{13}\beta}$ \\[0.4cm]

%$10^{-4}$  & \large $\alpha={\frac{7}{36} + \frac{7}{12}\beta}$ \\[0.4cm]

%$10^{-5}$  & \large $\alpha={136 -30\beta}$ \\[0.4cm]

$10^{-6}$  & \large $\alpha={\frac{56}{701} + \frac{427}{701}\beta}$ \\[0.4cm]

$10^{-7}$  & \large $\alpha={\frac{565}{289} + \frac{54}{289}\beta}$ \\[0.4cm]

$10^{-8}$  & \Large $\alpha=\frac{\beta}{    \log({\sfrac{59}{11} - \sfrac{13}{132}\beta})}$ \\[0.4cm]

$10^{-9}$  & \Large $\alpha=\frac{ 3^{{\sfrac{36}{73}}} 5^{\sfrac{2}{73}} 7^{\sfrac{25}{73}} \beta ^{{\sfrac{33}{73}}}} {2^{\sfrac{95}{73}}}$ \\[0.4cm]

$10^{-10}$  & \Large $\alpha=\frac{{2^{\sfrac{3}{4}}3^{\sfrac{11}{18}}\beta^{\sfrac{47}{36}}}}{{5^{\sfrac{11}{36}}7^{\sfrac{5}{6}}}}$ \\[0.4cm]

$10^{-11}$  & \Large $\alpha=\frac{{2^{\sfrac{108}{91}}5^{\sfrac{48}{91}}7^{\sfrac{9}{13}}\beta^{\sfrac{5}{91}}}}{{3^{\sfrac{172}{91}}}}$ \\[0.4cm]

$10^{-12}$  &  $None$ \\

	\end{tabular}
	\caption{Algebraic relations between $\alpha=\Jimm(\sqrt[3]{2})$ and $\beta=\Jimm(\sqrt[3]{4})$ up to various error levels}
\label{pslq1}
\end{table}

It can be seen on \autoref{pslq1} that as the allowed error tolerance is reduced, the simplest expression relating $\Jimm(\sqrt[3]{2})$ and $\beta=\Jimm(\sqrt[3]{4})$ gets more complicated, up to the point where for error tolerance of $10^{-12}$ no relations can be found within given coefficient size constraint.

\paragraph{Experiment 2.}  As a second investigation, we search whether the relation between $\sqrt[3]{2}$ and $2 \times  {\sqrt[3]{2}}$ is mapped to algebraically related numbers under $\Jimm$ transformation.

\begin{table}[h!]
\center
	\begin{tabular}{ c r }
\hline
	Tolerance & Relation Expression  \\ \hline \\

%$10^{-2}$  & \large $\alpha=4 - \beta$ \\[0.4cm]
%$10^{-3}$  & \large $\alpha=4 - \beta$ \\[0.4cm]
%$10^{-4}$  & \large $\alpha=\frac{83}{25} + \frac{-11}{25}\beta$ \\[0.4cm]
%$10^{-5}$  & \large $\alpha=\frac{-14}{19} + \frac{55}{19}\beta$ \\[0.4cm]
%$10^{-6}$  & \large $\alpha=\frac{101}{135} + \frac{226}{135}\beta$ \\[0.4cm]
$10^{-7}$  & \Large $\alpha=\frac{\sqrt{\sfrac{553+\sqrt{234969}}{110}} }{ \sqrt{\beta}}$ \\[0.4cm]

$10^{-8}$  & \Large $\alpha=\frac{ 2^{\sfrac{2}{7}} 3^{\sfrac{29}{14}} \beta^{\sfrac{10}{7}}  }{  5^{\sfrac{9}{14}} 7^{\sfrac{5}{14}}   }$ \\[0.4cm]

$10^{-9}$  & \Large $\alpha=\frac{ 2^{\sfrac{4}{9}} 3^{\sfrac{32}{9}}  }{  5^{\sfrac{5}{9}} 7^{\sfrac{8}{9}} \beta^{\sfrac{26}{9}}   }$ \\[0.4cm]

$10^{-10}$  & \Large $\alpha=\frac{ 2^{\sfrac{49}{38}} 3^{\sfrac{63}{19}}  }{  5^{\sfrac{15}{19}} 7^{\sfrac{33}{38}}  \beta^{\sfrac{107}{38}}   }$ \\[0.4cm]

$10^{-11}$  & \Large $\alpha=\frac{ 7^{\sfrac{35}{39}} \beta^{\sfrac{12}{13}}  }{  2^{\sfrac{43}{39}} 3^{\sfrac{2}{39}} 5^{\sfrac{2}{39}}   }$ \\[0.4cm]

$10^{-12}$  & \Large $\alpha=\frac{ 5^{\sfrac{29}{37}} 7^{\sfrac{20}{111}} \beta^{\sfrac{140}{111}}  }{  2^{\sfrac{64}{111}} 3^{\sfrac{44}{111}}   }$ \\[0.4cm]

$10^{-13}$  & \Large $\alpha=\frac{ 2^{\sfrac{199}{95}} 7^{\sfrac{99}{95}} \beta^{\sfrac{10}{19}}  }{  3^{\sfrac{107}{95}} 5^{\sfrac{78}{95}}   }$ \\[0.4cm]

$10^{-14}$  & \Large $\alpha=\frac{ 2^{\sfrac{56}{405}} 3^{\sfrac{4}{405}} 7^{\sfrac{233}{405}} \beta^{\sfrac{379}{405}}  }{  5^{\sfrac{97}{405}}   }$ \\[0.4cm]

$10^{-15}$  &  $None$ \\
	\end{tabular}
	\caption{Algebraic relations between $\alpha=\Jimm(^3\sqrt{2})$ and $\beta=\Jimm(2\times  {^3\sqrt{2}})$ up to various error levels}
\label{pslq2}
\end{table}

\autoref{pslq2} displays results similar to those of \autoref{pslq1}, which allows us to conclude that even if there still is a relation between $\Jimm$ images of those two numbers, the expression for it would be non-trivial.

\paragraph{Experiment 3.}  Finally, we search for a possible relation between the triple of numbers consisting of $(\alpha=\Jimm(^3\sqrt{2}),\beta=\Jimm(2\times  {^3\sqrt{2}}),\gamma=\Jimm(\frac{ {^3\sqrt{2}}}{2}))$, this time by building the dictionary with transformations and pairwise combinations of those three.

\begin{table}[h!]
\center
	\begin{tabular}{ c r }
\hline
	Tolerance & Relation Expression  \\ \hline \\

%$10^{-2}$  & \large $\alpha=\frac{1+\sqrt{21}}{2}$\\[0.4cm]
%$10^{-3}$  & \large $\alpha=\frac{8+\sqrt{960}}{14}$\\[0.4cm]
%$10^{-4}$  & \large $\alpha=\frac{22-\sqrt{28}}{6}$\\[0.4cm]
%$10^{-5}$  & \large $\alpha=\frac{1}{2} + \frac{13}{5}\beta - \frac{39}{10}\gamma$\\[0.4cm]
%$10^{-6}$  & \large $\alpha=\frac{7}{11} + \frac{37}{22}\beta + \frac{5}{11}\gamma$\\[0.4cm]
$10^{-7}$  & \large $\alpha=\frac{128}{227} + \frac{392}{227}\beta + \frac{121}{227}\gamma$\\[0.4cm]
$10^{-8}$  & \large $\alpha=\frac{117}{74} + \frac{153}{74}\beta - \frac{431}{74}\gamma$\\[0.4cm]
$10^{-9}$  & \large $\alpha=\frac{-64}{77} + \frac{694}{231}\beta - \frac{40}{231}\gamma$\\[0.4cm]
$10^{-10}$  & \large $\alpha=\frac{-88}{433} + \frac{908}{433}\beta + \frac{840}{433}\gamma$\\[0.4cm]
$10^{-11}$  & \Large $\alpha=\frac{\gamma}{ \sfrac{-88}{243} + \sfrac{61}{243}\beta + \sfrac{11}{18}\gamma }$\\[0.4cm]
$10^{-12}$  & \Large $\alpha=\frac{7^{\sfrac{31}{5}}  \gamma^{\sfrac{4}{5}}  }{  2^{\sfrac{7}{5}} 3^{\sfrac{7}{5}} 5^4 \beta^{\sfrac{23}{5}}  }$\\[0.4cm]
$10^{-13}$  &  $None$\\
	\end{tabular}
	\caption{Algebraic relations between $\alpha=\Jimm(^3\sqrt{2}),\beta=\Jimm(2\times  {^3\sqrt{2}}),$ and $\gamma=\Jimm(\frac{ {^3\sqrt{2}}}{2})$ up to various error levels}
\label{pslq3}
\end{table}

Once again we see on \autoref{pslq3} that the complexity of a possible relation between the triple of numbers grows with the allowed error tolerance up to the point where the space of expressions defined by constraints is exhausted.

\bigskip\noindent
{\bf Acknowledgements.} 
This research is sponsored by the T\"{U}B\.{I}TAK grant 115F412 and by the Galatasaray University research grant 17.504.001. 

%\nt{buraya birkaç modern referans daha eklersek iyi olur. Bir de referanslari duzneleylim.}
%
%\begin{thebibliography}{99}
%
%\bibitem{adam1}
%Continued fractions and transcendental numbers
%Boris Adamczewski, Yann Bugeaud, and Les Davison
%
%\bibitem{bugeaud}
%Exponents of Diophantine approximation
%Yann Bugeaud
%
%\bibitem{khinchinebook}
%Khintchine, Aleksandr Yakovlevich. Continued fractions. 1963.
%
%\bibitem{sympy}
%\url{http://sympy.org}
%
%\bibitem{mpmath}
%\url{http://mpmath.org/}
%
%\bibitem{fergy}
%Rabago, Julius Fergy T. "On k-Fibonacci Numbers with Applications to Continued Fractions." Journal of Physics: Conference Series. Vol. 693. No. 1. IOP Publishing, 2016.
%
%\bibitem{mpmathidentify}
%\url{http://docs.sympy.org/dev/modules/mpmath/identification.html}
%
%\bibitem{hakami}
%Hakami, Ali H. "An Application of Fibonacci Sequence on Continued Fractions." International Mathematical Forum. Vol. 10. No. 2. 2015.
%
%\bibitem{lebesgue}
%Uluda\u g, A. Muhammed, and Hakan Ayral. "A subtle symmetry of Lebesgue's measure" arXiv preprint arXiv:1501.03787.
%
%\bibitem{overlooked}
%Uluda\u g, A. Muhammed, and Hakan Ayral. "Jimm, a fundamental involution. 2015." arXiv preprint arXiv:1501.03787.
%\end{thebibliography}

\clearpage

\bibliographystyle{plain}
\bibliography{experimentalReferences}

\paragraph{Appendix-I: Maple code to evaluate Jimm on $\Q$.}
\begin{verbatim}
>with(numtheory)
>jimm := proc (q) local M, T, U, i, x; 
	T := matrix([[1, 1], [1, 0]]); 
	U := matrix([[0, 1], [1, 0]]); 
	M := matrix([[1, 0], [0, 1]]); 
	x := cfrac(q, quotients); 
	if x[1] = 0 then for i from 2 to nops(x) do 
		M := evalm(`&*`(`&*`(M, T^x[i]), U)) end do; 
	return M[2, 2]/M[1, 2] else for i to nops(x) do 
		M := evalm(`&*`(`&*`(M, T^x[i]), U)) end do; 
	return M[1, 2]/M[2, 2] end if 
end proc;
\end{verbatim}

\noindent
\paragraph{Appendix-II} $\Jimm$-transforms of some quadratic surds
$$
\begin{array}{|l|l|}
\hline N & \Jimm(\sqrt{N})\\
\hline 3& {1\over2}(\sqrt{13}+3)\\
\hline 5& {1\over3}(\sqrt{10}+1)\\
\hline 6& {1\over14}(\sqrt{221}+5)\\
\hline 7& {1\over6}(\sqrt{37}+1)\\
\hline 8& {1\over4}(\sqrt{17}+1)\\
\hline 10& {1\over7}(\sqrt{65}+4)\\
\hline 11& {1\over26}(\sqrt{901}+15)\\
\hline 12& {1\over34}(\sqrt{1517}+19)\\
\hline 13& {{1\over3}(\sqrt{13}}+2)\\
\hline 14& {1\over5}(\sqrt{34}+3)\\
\hline 15& {1\over18}(\sqrt{445}+11)\\
\hline 17& {1\over19}(\sqrt{442}+9)\\
\hline 18& {1\over78}(\sqrt{7453}+37)\\
\hline 19& {1\over730}(\sqrt{656101}+351)\\
\hline 20& {1\over23}(5\sqrt{26}+11)\\
\hline
\end{array}
\begin{array}{|l|l|}
\hline N & \Jimm(\sqrt{N})\\
\hline 21& {1\over307}(\sqrt{113570}+139)\\
\hline 22& {13\over307}(\sqrt{677}+142)\\
\hline 23& {1\over24}(\sqrt{697}+11)\\
\hline 24& {1\over50}(\sqrt{3029}+23)\\
\hline 26& {1\over49}(\sqrt{3026}+25)\\
\hline 27& {1\over194}(\sqrt{47437}+99)\\
\hline 28& {1\over139}(\sqrt{24362}+71)\\
\hline 29& {1\over495}(\sqrt{308026}+251)\\
\hline 30& {1\over238}(\sqrt{71285}+121)\\
\hline 31& {1\over17226}(\sqrt{376748101}+8945)\\
\hline 32& {1\over94}(\sqrt{11237}+49)\\
\hline 33& {1\over101}(\sqrt{12905}+52)\\
\hline 34& {1\over130}(\sqrt{21389}+67)\\
\hline 35& {1\over64}(\sqrt{5185}+33)\\
\hline 37& {1\over129}(\sqrt{20737}+64)\\
%\hline 38& {17\over518}(\sqrt{1157}+257)\\
%\hline 39& {1\over530}(\sqrt{350069}+263)\\
\hline
\end{array}
$$
\end{document}